\documentclass{amsart}
\usepackage{amsmath,amssymb,amsthm,amsfonts,amscd}
\usepackage[numeric]{amsrefs}

\newcommand{\nc}{\newcommand}
\newcommand{\rc}{\renewcommand}
\nc{\cmt}[1]{\noindent\textRed{{\fbox{W}} #1}\textBlack}
\nc{\cmtk}[1]{\noindent\textBlue{{\fbox{K}} #1}\textBlack}
\nc{\cmtm}[1]{\noindent\textGreen{{\fbox{M}} #1}\textBlack}

\nc{\on}{\operatorname}

\nc{\al}{{\alpha }}
\nc{\be}{{\beta }}
\nc{\ga}{{\gamma }}
\nc{\de}{{\delta }}
\nc{\del}{{\partial }}
\nc{\ep}{{\varepsilon }}
\nc{\vap}{{\epsilon }}

\nc{\ze}{{\zeta }}
\nc{\et}{{\eta }}
\rc{\th}{{\theta }}

\nc{\vth}{{\vartheta }}

\nc{\io}{{\iota }}
\nc{\ka}{{\kappa }}
\nc{\la}{{\lambda }}
\nc{\vrho}{{\varrho}}
\nc{\si}{{\sigma }}
\nc{\ups}{{\upsilon }}
\nc{\vphi}{{\varphi }}
\nc{\om}{{\omega }}

\nc{\Ga}{{\Gamma }}
\nc{\De}{{\Delta }}
\nc{\nab}{{\nabla}}
\nc{\Th}{{\Theta }}
\nc{\La}{{\Lambda }}
\nc{\Si}{{\Sigma }}
\nc{\Ups}{{\Upsilon }}
\nc{\Om}{{\Omega }}

\renewcommand\a{\alpha}

\renewcommand\l{\lambda}

\renewcommand\om{\omega}

\nc{\fA}{{\mathfrak A}}
\nc{\fB}{{\mathfrak B}}
\nc{\fC}{{\mathfrak C}}
\nc{\fD}{{\mathfrak D}}
\nc{\fE}{{\mathfrak E}}
\nc{\fF}{{\mathfrak F}}
\nc{\fG}{{\mathfrak G}}
\nc{\fH}{{\mathfrak H}}
\nc{\fI}{{\mathfrak I}}
\nc{\fJ}{{\mathfrak J}}
\nc{\fK}{{\mathfrak K}}
\nc{\fL}{{\mathfrak L}}
\nc{\fM}{{\mathfrak M}}
\nc{\fN}{{\mathfrak N}}
\nc{\fO}{{\mathfrak O}}
\nc{\fP}{{\mathfrak P}}
\nc{\fQ}{{\mathfrak Q}}
\nc{\fR}{{\mathfrak R}}
\nc{\fS}{{\mathfrak S}}
\nc{\fT}{{\mathfrak T}}
\nc{\fU}{{\mathfrak U}}
\nc{\fV}{{\mathfrak V}}
\nc{\fW}{{\mathfrak W}}
\nc{\fZ}{{\mathfrak Z}}
\nc{\fX}{{\mathfrak X}}
\nc{\fY}{{\mathfrak Y}}
\nc{\fa}{{\mathfrak a}}
\nc{\fb}{{\mathfrak b}}
\nc{\fc}{{\mathfrak c}}
\nc{\fd}{{\mathfrak d}}
\nc{\fe}{{\mathfrak e}}
\nc{\ff}{{\mathfrak f}}
\nc{\fg}{{\mathfrak g}}
\nc{\fh}{{\mathfrak h}}
\nc{\fiI}{{\mathfrak i}}  
\nc{\ffi}{{\mathfrak i}}  
\nc{\fj}{{\mathfrak j}}
\nc{\fk}{{\mathfrak k}}
\nc{\fl}{{\mathfrak{l}}}
\nc{\fm}{{\mathfrak m}}
\nc{\fn}{{\mathfrak n}}
\nc{\fo}{{\mathfrak o}}
\nc{\fp}{{\mathfrak p}}
\nc{\fq}{{\mathfrak q}}
\nc{\fr}{{\mathfrak r}}
\nc{\fs}{{\mathfrak s}}
\nc{\ft}{{\mathfrak t}}
\nc{\fu}{{\mathfrak u}}
\nc{\fv}{{\mathfrak v}}
\nc{\fw}{{\mathfrak w}}
\nc{\fz}{{\mathfrak z}}
\nc{\fx}{{\mathfrak x}}
\nc{\fy}{{\mathfrak y}}

\nc{\bA}{{\mathbb A}}
\nc{\bB}{{\mathbb B}}
\nc{\bC}{{\mathbb C}}
\nc{\bD}{{\mathbb D}}
\nc{\bE}{{\mathbb E}}
\nc{\bF}{{\mathbb F}}
\nc{\bG}{{\mathbb G}}
\nc{\bH}{{\mathbb H}}
\nc{\bI}{{\mathbb I}}
\nc{\bJ}{{\mathbb J}}
\nc{\bK}{{\mathbb K}}
\nc{\bL}{{\mathbb L}}
\nc{\bM}{{\mathbb M}}
\nc{\bN}{{\mathbb N}}
\nc{\bO}{{\mathbb O}}
\nc{\bP}{{\mathbb P}}
\nc{\bQ}{{\mathbb Q}}
\nc{\bR}{{\mathbb R}}
\nc{\bS}{{\mathbb S}}
\nc{\bT}{{\mathbb T}}
\nc{\bU}{{\mathbb U}}
\nc{\bV}{{\mathbb V}}
\nc{\bW}{{\mathbb W}}
\nc{\bZ}{{\mathbb Z}}
\nc{\bX}{{\mathbb X}}
\nc{\bY}{{\mathbb Y}}

\nc{\cA}{{\mathcal A}}
\nc{\cB}{{\mathcal B}}
\nc{\cC}{{\mathcal C}}
\nc{\cD}{{\mathcal D}}
\nc{\cE}{{\mathcal E}}
\nc{\cF}{{\mathcal F}}
\nc{\cH}{{\mathcal H}}
\nc{\cI}{{\mathcal I}}
\nc{\cJ}{{\mathcal J}}
\nc{\cK}{{\mathcal K}}
\nc{\cL}{{\mathcal L}}
\nc{\cM}{{\mathcal M}}
\nc{\cN}{{\mathcal N}}
\nc{\cO}{{\mathcal O}}
\nc{\cP}{{\mathcal P}}
\nc{\cQ}{{\mathcal Q}}
\nc{\cR}{{\mathcal R}}
\nc{\cS}{{\mathcal S}}
\nc{\cT}{{\mathcal T}}
\nc{\cU}{{\mathcal U}}
\nc{\cV}{{\mathcal V}}
\nc{\cW}{{\mathcal W}}
\nc{\cZ}{{\mathcal Z}}
\nc{\cX}{{\mathcal X}}
\nc{\cY}{{\mathcal Y}}

\newcommand\codim{\operatorname{codim}}

\newcommand{\gr}{\operatorname{gr}}

\newcommand{\gobble}[1]{}
  \newcommand{\rangeref}[2]{%
    \ref{#1}--\afterassignment\gobble\fam 0\ref{#2}
  }

\begin{document}

\title{Hodge theory and unitary representations, in the example of $SL(2,\bR)$.}

\author{Wilfried Schmid}
\address{Department of Mathematics, Harvard University, Cambridge, MA 02138}
\email{schmid@math.harvard.edu}
\thanks{The first author was supported in part by NSF grant DMS-1300185.}

\author{Kari Vilonen}\address{Department of Mathematics, Northwestern University, Evanston, IL 60208, also Department of Mathematics, Helsinki University, Helsinki, Finland}
\email{vilonen@northwestern.edu, kari.vilonen@helsinki.fi}
\thanks{The second author was supported in part by NSF grants DMS-1402928, DMS-1069316, and the Academy of Finland.}



\dedicatory{Dedicated to David Vogan, on the occasion of his sixtieth birthday}

\maketitle

In our paper \cite{vs:2011} we formulated a conjecture on unitary representations of reductive Lie groups. We are currently working towards a proof; the technical difficulties are formidable. It has been suggested that an explicit description in the case of $SL(2,\bR)$ would be helpful. The unitary representations of $SL(2,\bR)$ have been worked out in great detail, of course, but even in this special case our construction of the inner product in terms of the $\cD$-module realization is not obvious.

We begin with a quick summary of our conjecture in the general case of a reductive, linear, connected Lie group $G_\bR$, with maximal compact subgroup $K_\bR$. We let $G$ and $K$ denote the complexifications. The complex group $G$ contains a unique compact real form $U_\bR$ such that $U_\bR \cap G_\bR = K_\bR$. Then $U_\bR$ acts transitively on the flag variety $X$ of $G$, and $K$ acts with finitely many orbits. The points of $X$ corres\-pond to Borel subalgebras $\fb$ of the Lie algebra $\fg$ of $G$. The quotients $\fh=\fb/[\fb,\fb]$ constitute the fibers of a flat vector bundle. Since $X$ is simply connected, we can think of $\fh$ as a fixed vector space. This is the ``universal Cartan", and is acted on by the ``universal Weyl group" $W$. Its dual $\fh^*$ contains the ``universal root system" $\Phi$ and the system of positive roots $\Phi^+$, chosen so that $[\fb,\fb]$ becomes the direct sum of the negative root spaces; $\fh^*$ also contains the ``universal weight lattice" $\La$. Further notation: lower case Fraktur letters, such as $\fg_\bR$, $\fk_\bR$, $\fg$, refer to the Lie algebras of $G_\bR$, $K_\bR$, $G$, etc.

Via the Harish Chandra isomorphism, $\fh^*$ parameterizes the characters $\chi_\lambda$ of the center of the universal enveloping algebra, with $\chi_\l=\chi_\mu$ if and only if $\mu=w\l$ for some $w\in W$. We shall say that a Harish Chandra module $M$ has a real infinitesimal character if it is of the form $\chi_\la$ with $\la\in \bR\otimes_\bZ \La$. David Vogan, many years ago, pointed out that to understand the irreducible unitary representations of $G_\bR$ it suffices to treat the case of real infinitesimal character \cite{knapp:1986}. Let then $M_\la$ be an irreducible Harish Chandra module with real infinitesimal character $\chi_\la$. Since $\la$ is determined only up to the Weil group action, we may and shall assume that $\la$ is dominant, i.e.,
\begin{equation}
\label{dominant}
(\a,\la)\ \geq \ 0\ \ \ \text{for all}\ \ \a\in\Phi^+.
\end{equation}
To determine whether or not $M_\la$ underlies an irreducible unitary representation, one needs to know if it carries a nonzero $\fg_\bR$-invariant hermitian form $(\ ,\ )_{\fg_\bR}$ -- this question has a simple answer, see below -- and, when that is the case, if $(\ ,\ )_{\fg_\bR}$ has a definite sign.

Vogan and his coworkers \cite{vogan:2009} made the important observation that the condition of having a real infinitesimal character ensures the existence of a nonzero $\fu_\bR$-invariant hermitian form $(\ ,\ )_{\fu_\bR}$. If both types of hermitian forms exist, they are explicitly related: the Cartan involution $\,\th: \fg\to\fg\,$ then acts also on the Harish Chandra module $M_\la$, and after suitable rescaling of the hermitian forms,
\begin{equation}
\label{gvsu}
(v_1,v_2)_{\fg_\bR}\ =\ (\th v_1,v_2)_{\fg_\bR}\ .
\end{equation}
The $\fu_\bR$-invariant form is easier to deal with, both computationally and from a geometric point of view. At the same time the action of $\,\th\,$ on $M_\la$ can be described quite concretely. Thus, if one understands the hermitian form $(\ ,\ )_{\fu_\bR}$, one can decide if $M_\la$ is unitarizable.

As usual we write $\rho$ for the half sum of the positive roots. We let $\cD$ denote the sheaf of linear differential operators, with algebraic coefficients, on the flag variety $X$; here $X$ is equipped with the Zariski topology. The sheaf of algebras $\cD$ can be twisted by $G$-equivariant line bundles, and more generally, by any $\la\in\fh^*$. It is convenient to parameterize the twists so that $\cD_\la$, for $\la\in\La+\rho$, acts on sections of the $G$-equivariant line bundle $\cL_{\la-\rho}\to X$ with Chern class $\la-\rho\in\La\cong H^2(X,\bZ)$; for arbitrary $\la\in\fh$ one then defines $\cD_\la$ by a process of analytic continuation. The sheaves $\cD_\la$ are $G$-equivariant, in the Zariski sense locally isomorphic to $\mathcal D$, and every $\zeta\in \fg$ acts as an infinitesimal automorphism and thus defines a global section of $\cD_\la$. Note that $\cD_\rho=\cD$, and that $\cD_{-\rho}$ acts on sections $\cL_{-2\rho}\,$= canonical bundle of $X$.

The Beilinson-Bernstein construction\begin{footnote}{A more detailed summary of the Beilinson-Bernstain construction of Harish Chandra modules can be found in \cite{HMSW:1987}}\end{footnote} realizes the irreducible Harish Chandra module $M_\la$, with $\la$ real and dominant as in (\ref{dominant}), as the space of global sections
\begin{equation}
\label{BB1}
M_\la \ = \ H^0(X,\cM_\la)
\end{equation}
of an irreducible, $K$-equivariant sheaf of $\cD_\la$-modules $\cM_\la$ -- or for short, an irreducible, $K$-equivariant $\cD_\la$-module. Then $\fg$ acts on $M_\la \ = \ H^0(X,\cM_\la)$ via the inclusion $\fg\hookrightarrow H^0(X,\cD_\la)$.  The correspondence between the Harish Chandra module $M_\la$ and the ``Harish Chandra sheaf" $\cM_\la$ extends functorially to all Harish Chandra modules of finite length, with infinitesimal cha\-racter $\chi_\la$. Irreducible Harish Chandra sheaves $\cM_\la$ are easy to describe: they arise from direct images, in the category of $\cD_\la$-modules, under the embedding $j:Q\hookrightarrow X$ of a $K$-orbit $Q$ in $X$, applied to a $K$-equivariant ``twisted local system" $\bC_{Q,\la}$ on $Q$, with twist $\la-\rho$. A formal, general definition of $\bC_{Q,\la}$ would lead too far; in the special case of $G_\bR= SL(2,\bR)$ we describe it implicitly in (\ref{C*construction1}) below, where its generating section will be denoted by $\sigma_0^{\frac{\la - 1}{2}}$. In any case, the tensor product $\cO_Q\otimes_\bC\bC_{Q,\la}$ has the structure of a $\cD_{Q,\la}$-module on the $K$-orbit $Q$, and the direct image $j_*(\cO_Q\otimes_\bC \bC_{Q,\la})$ that of a Harish Chandra sheaf: a $K$-equivariant $\cD_\la$-module on $X$. In general the direct image is not irreducible, but it contains a unique irreducible $\cD_\la$-submodule\begin{footnote}{Since we assumed $G_\bR$, and hence also $K$, to be connected, any $\cD_\la$-subsheaf of a Harish Chandra sheaf is automatically $K$-equivariant and is therefore also a Harish Chandra sheaf.}\end{footnote}, and
\begin{equation}
\label{BB2}
\cM_\la \ = \ \text{unique irreducible submodule of}\,\ j_*(\cO_Q\otimes_\bC\bC_{Q,\la})\ .
\end{equation}
The realization (\ref{BB2}) of irreducible Harish Chandra sheaves is unique. It almost sets up a bijection between irreducible Harish Chandra modules $M_\la$, with the parameter $\la$ of the infinitesimal character as in (\ref{dominant}), and $K$-equivariant, irreducible local systems $\bC_{Q,\la}$, with twist $\la-\rho$, on $K$-orbits $Q\subset X$ -- the qualifier ``almost" is necessary because when $\la$ is singular, certain irreducible Harish Chandra sheaves have no nonzero global sections. This phenomenon explains why the classification of irreducible Harish Chandra modules with {\it regular infinitesimal character} looks simpler than that of irreducible Harish Chandra modules with {\it singular infinitesimal character}.

We shall not attempt to summarize Saito's theory of mixed Hodge modules here. Rather, we shall state the relevant facts, which apply to all members of the category of ``geometrically constructible" Harish Chandra sheaves $\cM_\la$ -- this includes in particular the sheaves obtained by the standard $\cD$-module operations applied to $\cD_\la$-modules of the type $j_*(\cO_Q\otimes_\bC\bC_{Q,\la})$ and their $\cD_\la$-subsheaves. A mild generalization of Saito's theory\begin{footnote}{Without the assumption of an underlying rational structure, which Saito requires.}\end{footnote} puts three additional structures on each object $\cM_\la$\,. First of all, the {\it weight filtration}, a functorial, finite increasing filtration
\begin{equation}
\label{wf1}
0\ \subset W_0\cM_\la \ \subset \  W_1\cM_\la \ \subset \ \cdots\ \subset \ W_k\cM_\la \ \subset \  \cdots \ \subset \ W_n\cM_\la \ = \ \cM_\la
\end{equation}
by $\cD_\la$-subsheaves, with completely reducible quotients $\,W_k\cM_\la/W_{k-1}\cM_\la\,$ which are themselves objects in the category of Harish Chandra sheaves. Secondly, the {\it Hodge filtration}, a typically infinite, increasing filtration
\begin{equation}
\label{hf1}
0\, \subset F_a\cM_\la \, \subset \, \cdots\, \subset \, F_p\cM_\la \, \subset \, F_{p+1}\cM_\la\, \subset \, \cdots \, \subset \, \cM_\la\, = \, {\cup}_{p\geq a}\,F_p\cM_\la
\end{equation}
by $\cO_X$-coherent, $K$-equivariant, $\cO_X$-submodules. This is a good filtration in the sense of $\cD$-module theory: let $(\cD_\la)_d\subset \cD_\la$ denote the $\cO_X$-subsheaf of differential operators of degree at most $d$; then
\begin{equation}
\label{gf}
(\cD_\la)_d\,F_p\cM_\la\ \subseteq\ F_{p+d}\cM_\la\,,\ \ \ \text{with equality holding if}\ \ p\gg 0\,.
\end{equation}
The third ingredient, the {\it polarization} on any irreducible Harish Chandra sheaf $\cM_\la$, is a nontrivial $\cD_\la\times \overline \cD_\la$-bilinear pairing
\begin{equation}
\label{pol1}
P\,:\,\cM_\la \times \overline \cM_\la\ \ \longrightarrow\ \ \cC^{-\infty}(X_\bR)\,.
\end{equation}
Here $\cC^{-\infty}(X_\bR)$ refers to the sheaf of distributions on $X$, considered as $C^\infty$ manifold, and $\overline \cM_\la$ is the complex conjugate of $\cM_\la$, viewed as $\overline \cD_\la$-module on $X$, equipped with the complex conjugate algebraic structure.

Morphisms in the category of mixed Hodge modules preserve both filtrations strictly: if $T:\cM \to \cN$ is a morphism, then $T(F_p\cM) = (T\,\cM)\cap(F_p\cN)$, and ana\-logously for the weight filtration. We should also mention Saito's normalization of the indexing of the two filtrations. Going back to (\ref{BB2}), the Hodge filtration on the sheaf $\cO_Q\otimes_\bC\bC_{Q,\la}$ on $Q$ is trivial, in the sense that $F_0(\cO_Q\otimes_\bC\bC_{Q,\la})=\cO_Q\otimes_\bC\bC_{Q,\la}$ and $F_{-1}(\cO_Q\otimes_\bC\bC_{Q,\la})=0$. As a sheaf on $Q$ it is irreducible and has weight equal to $\dim Q$. The process of direct image shifts the lowest index of the Hodge filtration to $a=\codim Q$, and puts the weights into degrees $\geq \dim Q$, with the irreducible subsheaf $\cM_\la$ having weight equal to $\dim Q$.

The polarization leads to a geometric description of the $\fu_\bR$-invariant hermitian form on any irreducible Harish Chandra module $M_\la$. Let $\om$ denote the -- unique, up to scaling -- $U_\bR$-invariant measure on $X_\bR$. Like any smooth measure on a compact $C^\infty$ manifold it can be integrated against any distribution. When $M_\la$ is realized as the space of global sections of the corresponding Harish Chandra sheaf $\cM_\la$ as in (\ref{BB1}), then
\begin{equation}
\label{pol2}
(s_1,s_2)_{\fu_\bR}\ = \ \int_{X_\bR} P(s_1,\bar s_2)\,\om\ \ \ \text{for $\,s_1,s_2\in H^0(X,\cM_\la)$}\,,
\end{equation}
does indeed define a $\fu_\bR$-invariant hermitian form, as can be checked readily \cite{vs:2011}. The Cartan involution $\,\theta\,$ acts on $X$ and on the set of $K$-orbits in $X$. If $\,\theta\,$ fixes a particular $K$-orbit $Q$, then it acts on the twisted local systems on $Q$, and if it fixes also the twisted local system $\bC_{Q,\la}$ as in (\ref{BB2}), then it acts on the sections of the direct image $j_*(\cO_Q\otimes_\bC\bC_{Q,\la})$ and of its unique irreducible subsheaf $\cM_\la$. In this sense, the action of $\,\theta\,$ on the Harish Chandra module $M_\la$ -- which relates $(\ ,\ )_{\fu_\bR}$ to $(\ ,\ )_{\fg_\bR}$ as in (\ref{gvsu}) -- is visible geometrically.

Via the global section functor the Hodge and weight filtrations induce filtrations on $M_\la$, the space of global sections of the Harish Chandra sheaf $\cM_\la$, wether or not the latter is irreducible:
\begin{equation}
\label{wf2hf2}
\begin{aligned}
&0\ \subset W_0 M_\la \ \subset \  W_1 M_\la \ \subset \ \cdots\ \subset \ W_k M_\la \ \subset \  \cdots \ \subset \ W_n M_\la \ = \ M_\la\,,
\\
&0\, \subset F_a M_\la \, \subset \, \cdots\, \subset \, F_p M_\la \, \subset \, F_{p+1} M_\la\, \subset \, \cdots \, \subset \, M_\la\, = \, {\cup}_{p\geq a}\,F_p M_\la\,.
\end{aligned}
\end{equation}
The $W_kM_\la$ are Harish Chandra submodules of $M_\la$, and the $F_p M_\la$ are finite dimensional, $K$-invariant subspaces. In the irreducible case the weight filtration collapses as was mentioned earlier, $M_\la$ in (\ref{BB1},\ref{BB2}) has weight $\dim Q$, and the lowest index in the Hodge filtration is $a=\codim Q$. We can now state our conjecture. It asserts that if $M_\la$ is irreducible, the $\fu_\bR$-invariant hermitian form is nondegenerate on each $F_p M_\la$, and
\begin{equation}
\label{conj1}
(-1)^{p-a}\,(s,s)_{\fu_\bR}\ > \ 0 \ \ \text{for all nonzero $\,s\in F_p M_\la\cap (F_{p-1}M_\la)^\perp$}\,.
\end{equation}
Whenever $M_\la$ also admits a $\fg_\bR$-invariant hermitian form it would be related to $\fu_\bR$-invariant one via (\ref{gvsu}), and the resulting hermitian form $(\ ,\ )_{\fg_\bR}$ would then have a definite sign if and only if $M_\la$ is unitarizable.

The significance of the conjecture is discussed in \cite{vs:2011}. While it does not amount to a description of the unitary dual of $G_\bR$ in terms of representation parameters, it puts the study of the irreducible unitary representation into a functorial context.

We now turn to the example of $SL(2,\bR)$. It is conjugate to $SU(1,1)$ under an inner automorphism of $SL(2,\bC)$, and various formulas have a simpler appearance for $SU(1,1)$. Thus we suppose $G=SL(2,\bC)$,
\begin{equation}
\label{groups}
\begin{aligned}
&G_\bR\ = \ SU(1,1)\ = \ \left. \left\{
\begin{pmatrix} \a  &  \be \\ - \bar{\be} &  \bar{\a} \end{pmatrix}
\ \right| \ \a,\,\be\in \bC\,,\,\ |\a|^2-|\be|^2=1\  \right\},
\\
&K_\bR\ = \ \left. \left\{
\begin{pmatrix} \a  &  0 \\ 0 &  \bar{\a} \end{pmatrix}
\ \right| \ \a\in \bC\,,\,\ |\a|=1\  \right\},\, \ \ K\ = \ \left. \left\{
\begin{pmatrix} \a  &  0 \\ 0 &  \a^{-1} \end{pmatrix}
\ \right| \ \a\in \bC^*\  \right\},
\end{aligned}
\end{equation}
and $U_\bR = SU(2)$. These groups act on the flag variety of $G$,
\begin{equation}
\label{SL2flag}
X\ \ = \ \ \bP^1\ \ = \ \ \bC \cup \{\infty\}\,,
\end{equation}
by linear fractional transformations, and $K$ acts with three orbits, namely $\{0\}$, $\{\infty\}$, and $\bC^*$. In the notation of (\ref{groups},\ref{SL2flag}), $K$ acts on $\bC^*$ by $\a^2$, so $\bC^*$ admits two irreducible $K$-equivariant local systems, corresponding to the trivial and the nontrivial character of the (component group of the) generic isotropy group $\{\pm 1\}$. This is true both in the scalar -- i.e., non-twisted -- and twisted case. Since $K$ is connected, the two point orbits admit only the trivial irreducible $K$-equivariant local system. The Cartan involution,
\begin{equation}
\label{SL2theta}
\theta\ \ = \ \ \text{conjugation by}\ \ \ \begin{pmatrix} i  &  0 \\ 0 &  -i \end{pmatrix}\,,
\end{equation}
is inner, it preserves each of the three orbits and the $K$-equivariant local systems on them. Thus all irreducible Harish Chandra modules with real infinitesimal character admit both $\fu_\bR$- and $\fg_\bR$-invariant hermitian forms.

The dual $\fh^*$ of the universal Cartan can be identified with $\bC$ so that $\La\cong \bZ$, $\Phi\cong \{\pm 2\}$, and $\rho \cong 1$. With this identification an infinitesimal character $\chi_\la$ is real in the earlier sense if and only if $\la\in \bR$, and $\la\in\bR$ is dominant if and only if $\la\geq 0$. The standard $SL_2$-triple
\begin{equation}
\label{sl2triple1}
e_+\ = \ \begin{pmatrix} 0  &  1 \\ 0 &  0 \end{pmatrix}\,,\ \ \ e_-\ = \ \begin{pmatrix} 0  &  0 \\ 1 &  0 \end{pmatrix}\,,\ \ \ h\ = \ \begin{pmatrix} 1  &  0 \\ 0 &  -1 \end{pmatrix}
\end{equation}
spans $\fg$ over $\bC$, with $\fk$ spanned by $h$, and satisfies the conjugation relations $\bar e_+=e_-$\,, $\bar h = -h$. The elements of this triple operate on the sheaf of algebraic functions on $\bC \cup \{\infty\}$ by infinitesimal left translation. One computes readily that via this action,
\begin{equation}
\label{sl2triple2}
e_+\ \cong \ -\frac{d\ }{d z}\,, \ \ \ e_-\ \cong \ z^2\,\frac{d\ }{d z}\,,\ \ \ h\ \cong \ -2z\,\frac{d\ }{d z}\,.
\end{equation}
For example, $e_+$ acts on $f(z)$ by the derivative with respect to $t$, at the origin, of $f(\exp(-t\,e_+)z)=f(z-t)$, resulting in the formula $(e_+f)(z)= -z\frac{d f}{d z}(z)$; the other cases are treated similarly.

The $G$-equivariant line bundle $\cL_2$ coincides with the tangent bundle of $\bP^1$, so we can identify (\ref{sl2triple2}) with a basis of the space of global sections of $\cL_2$. However, for notational reasons, we choose the new symbols $\sigma_2$, $\sigma_0$, $\sigma_{-2}$, corresponding to $e_+$, $h$, $e_-$, in that order. Then $\sigma_2$ vanishes to second order at $\infty$, $\sigma_{-2}$ vanishes to second order at $0$, and $\sigma_0$ has first order zeros at both $0$ and $\infty$, and these are the only zeroes in each case. Moreover,
\begin{equation}
\label{sl2triple3}
\begin{aligned}
&e_+\,\sigma_2\ = \ 0\,,\ \ \ \ \,\,h\,\sigma_2\ = \ 2\sigma_2\,,\,\ \ \ \ e_-\,\sigma_2\ = \ -\sigma_0\ = \ {\textstyle\frac 12}\,z\, \sigma_2\,,
\\
&e_+\, \sigma_0\ \, = \,\ -2\sigma_2\,\  = \,\ -z^{-1}\sigma_0\,,\ \ \ \ h\,\sigma_0\ = \ 0\,,\ \ \ \ \ \ e_-\,\sigma_0\ = \ 2\,\sigma_{-2}\ = \ -z \sigma_0\,,
\\
&e_+\,\sigma_{-2}\ = \ \sigma_0\ = \ -{\textstyle\frac 12}\,z^{-1} \sigma_{-2}\,,\ \ \ \ h\,\sigma_{-2}\ = \ -2\sigma_{-2}\,,\ \ \ \ \ e_-\,\sigma_{-2}\ = \ 0\,,
\end{aligned}
\end{equation}
as can be read off from  (\ref{sl2triple2}). The $U_\bR$-invariant measure on $\bP^1$ is
\begin{equation}
\label{URmetric1}
\om\ \ = \ \ (\,1\,+\,|z|^2\,)^{-2}\,dz\,d\bar z\,.
\end{equation}
The coefficient of $dz\,d\bar z$ in this formula can be interpreted as the squared length of $\frac{d\ }{dz}$ with respect to the $U_\bR$-invariant hermitian metric, or equivalently, the squared length of $\sigma_2$ relative to the $U_\bR$-invariant hermitian metric on the line bundle $\cL_2$. Thus
\begin{equation}
\label{URmetric2}
\|\sigma_2\|\, =\, \frac{1}{1\,+\,|z|^2}\ \,,\ \ \ \|\sigma_0\|\, =\, \frac{2\,|z|}{1\,+\,|z|^2}\ \,,\ \ \ \|\sigma_{-2}\|\, =\, \frac{|z|^2}{1\,+\,|z|^2}
\end{equation}
describes the length, as measured by the $U_\bR$-invariant metric on $\cL_2$, of the three sections $\sigma_2$, $\sigma_0$, $\sigma_{-2}$.

As was mentioned already, there exist two irreducible $K$-invariant local systems, with twist $\la-\rho$, on the $K$-orbit $\bC^*$, corresponding to the trivial and the nontrivial character of the generic isotropy subgroup $\{\pm 1\}$ of $K$. The corresponding Harish Chandra sheaves can be realized as
\begin{equation}
\label{C*construction1}
\begin{aligned}
&\cM_{\bC^*\!,\,\la,\text{even}}\ = \ \{\,f(z)\,\sigma_0^{\frac{\la-1}{2}}\ |\ f\in \bC(z)\,\}\,,
\\
&\qquad\cM_{\bC^*\!,\,\la,\text{odd}}\ = \ \{\,f(z)\,z^{1/2}\,\sigma_0^{\frac{\la-1}{2}}\ |\ f\in \bC(z)\,\}\,.
\end{aligned}
\end{equation}
These are Zariski-locally defined algebraic functions, multiplied by the ``section" $\sigma_0^{\frac{\la-1}{2}}$ of the formal power $\cL_2^{\frac{\la-1}{2}}$, either on $\bC^*$ (in the even case), where the section is well defined, or its twofold cover (in the odd case). As such they are naturally $\cD_\la$-modules on $\bC^*$, and then, via the direct image functor corresponding to the open embedding $\bC^*\subset\bC \cup \{\infty\}$, on all of $\bC \cup \{\infty\}$. How $\cD_\la$ acts is not so relevant for us, but the action of $\fg \subset \Ga \cD_\la$ is. That action is given by the product rule, with $\fg$ acting on $f(z)$ or $z^{1/2}\,f(z)$ according to the formulas (\ref{sl2triple2}), and on the formal powers of $\sigma_0$ according to (\ref{sl2triple3}). Since $\sigma_0$ has first order zeroes at $0$ and $\infty$,
\begin{equation}
\label{C*construction2}
\begin{aligned}
&f(z)\,\sigma_0^{\frac{\la-1}{2}}\ \sim \ f(z)\,z^{\frac{\la-1}{2}}\ \ \ \text{near the origin, and}
\\
&\qquad f(z)\,\sigma_0^{\frac{\la-1}{2}}\ \sim \ f(z)\,z^{-\frac{\la-1}{2}}\ \ \ \text{near $\infty$}\,.
\end{aligned}
\end{equation}
In particular, $\cM_{\bC^*\!,\,\la,\text{even}}$ is reducible if and only if $\la$ is an odd integer, whereas $\cM_{\bC^*\!,\,\la,\text{odd}}$ reduces if and only if $\la$ is even.

We recall that the sheaves (\ref{C*construction1}), when restricted to $\bC^*$, are irreducible and have weight one, which is the dimension of $\bC^*$. That remains correct for these sheaves on all of $\bC \cup \{\infty\}$ when they are irreducible:
\begin{equation}
\label{C*construction3}
\begin{aligned}
&W_0\,\cM_{\bC^*\!,\,\la,\text{even}}\, = \, 0\ \ \ \text{and}\ \ \ W_1\,\cM_{\bC^*\!,\,\la,\text{even}}\, = \, \cM_{\bC^*\!,\,\la,\text{even}}\ \ \ \text{if}\, \ \la\notin 2\bZ+1\,,
\\
&\qquad W_0\,\cM_{\bC^*\!,\,\la,\text{odd}}\, = \, 0\ \ \ \text{and}\ \ \ W_1\,\cM_{\bC^*\!,\,\la,\text{odd}}\, = \, \cM_{\bC^*\!,\,\la,\text{odd}}\ \ \ \text{if}\, \ \la\notin 2\bZ\,.
\end{aligned}
\end{equation}
In the reducible case,
\begin{equation}
\label{C*construction4}
\begin{aligned}
&W_0\,\cM_{\bC^*\!,\,2m+1,\text{even}}\ = \ 0\,,\ \ \ \ W_1\,\cM_{\bC^*\!,\,2m+1,\text{even}}\ = \ \cO_{\bP^1}(\cL_{2m})\,,
\\
&\qquad\text{and}\ \ W_2\,\cM_{\bC^*\!,\,2m+1,\text{even}}\ = \ \cM_{\bC^*\!,\,2m+1,\text{even}}\,;
\\
&\qquad\qquad\qquad W_0\,\cM_{\bC^*\!,\,2m,\text{odd}}\ = \ 0\,,\ \ \ W_1\,\cM_{\bC^*\!,\,2m,\text{odd}}\ = \ \cO_{\bP^1}(\cL_{2m-1})\,,
\\
&\qquad\qquad\qquad\qquad\text{and}\ \ W_2\,\cM_{\bC^*\!,\,2m,\text{odd}}\ = \ \cM_{\bC^*\!,\,2m,\text{odd}}\,.
\end{aligned}
\end{equation}
To justify these descriptions of the weight filtrations one should notice that $\sigma_0^m$ can be viewed as a section of $\cL_2^m=\cL_{2m}$, and $z^{1/2}\sigma_0^{1/2}$ as a meromorphic section of $\cL_1$. The quotients $\gr_{W,2}\cM_{\bC^*\!,\,2m+1,\text{even}}$ and $\gr_{W,2}\cM_{\bC^*\!,\,2m,\text{odd}}$ are Harish Chandra sheaves supported on $\{0,\infty\}$. We shall discuss these later.

The Hodge filtration for the sheaves (\ref{C*construction1}) starts at level $a=0$, since that is the codimension. In general the Hodge filtration of the direct image under an open embedding is governed by the -- not necessarily integral -- order of poles. The case of $\bC^*\hookrightarrow \bC$, and analogously for $\bC^*\hookrightarrow \bC^*\cup\{\infty\}$, is especially simple: poles of order $\leq 1$ have Hodge level $0$, those of order $\leq 2$ have Hodge level $1$, and so forth. Thus, in view of (\ref{C*construction1},\ref{C*construction2}), for $n\in \bZ$ and $p\geq 0$,
\begin{equation}
\label{C*construction5}
\begin{aligned}
&z^n\,\sigma_0^{\frac{\la-1}{2}}\, \in \, F_p\,\cM_{\bC^*\!,\,\la,\text{even}}\ \ \ \Longleftrightarrow\ \ \ -{\textstyle\frac{\la+1}{2}} - p \leq n\leq {\textstyle\frac{\la+1}{2}} + p\,,
\\
&\qquad z^{n+1/2}\sigma_0^{\frac{\la-1}{2}}\, \in \, F_p\,\cM_{\bC^*\!,\,\la,\text{odd}}\ \ \ \Longleftrightarrow\ \ \ -{\textstyle\frac{\la+1}{2}} - p \leq n + {\textstyle\frac{1}{2}}\leq {\textstyle\frac{\la+1}{2}} + p\,,
\end{aligned}
\end{equation}
for all $n\in \bZ$ and $p\geq 0$. In the reducible case, there is a connection between the induced Hodge filtrations on the quotient sheaves and the intrinsic Hodge filtrations on the quotients; this, too, will be described later.

We now turn to the polarizations of the sheaves $\cM_{\bC^*\!,\,2n+1,\text{even}}$, $\cM_{\bC^*\!,\,2n+1,\text{odd}}$ and the resulting hermitian forms on their spaces of global sections. On the $K$-orbit $\bC^*$ these sheaves are always irreducible, and the only possible hermitian pairing of the type (\ref{pol1}) on $\bC^*$ is, up to scaling,
\begin{equation}
\label{pol3}
P\left(\,f(z)\,\sigma_0^{\frac{\la-1}{2}}\,,\ \overline{g(z)\,\sigma_0^{\frac{\la-1}{2}}}\,\right)\ \ = \ \ f(z)\,\overline{g(z)}\,\|\sigma_0\|^{\la - 1}\ ,
\end{equation}
which is a real analytic function, and thus distribution, on $\bC^*$. This is correct in both cases, if we take $f,g\in \bC(z)$ in the case of $\cM_{\bC^*\!,\,2n+1,\text{even}}$, and $f,g\in z^{1/2}\bC(z)$ in the case of $\cM_{\bC^*\!,\,2n+1,\text{odd}}$; (\ref{URmetric2}) makes the last factor on the right explicit. The two sheaves were defined as the $\cD_\la$-module direct image under the open embedding $\bC^*\hookrightarrow \bC \cup \{\infty\}$ which, as always in the case of open embeddings, coincides with the $\cO$-module direct image. The general theory ensures that
\begin{equation}
\label{pol4}
f(z)\,\overline{g(z)}\,\|\sigma_0\|^{\la - 1}
\end{equation}
makes sense as global distribution on $\bC\cup\{\infty\}$, for all global sections $f(z)\,\sigma_0^{\frac{\la-1}{2}}$, $g(z)\,\sigma_0^{\frac{\la-1}{2}}$, provided the sheaf in question is irreducible.

To see how this works out in the current setting, applied to the spaces of global sections $M_{\bC^*\!,\,2n+1,\text{even}}$, $M_{\bC^*\!,\,2n+1,\text{odd}}$ of the two sheaves, we note that
\begin{equation}
\label{basis1}
\begin{aligned}
&M_{\bC^*\!,\,\la,\text{even}}\ \ \ \ \text{has basis}\ \ \ \{\,z^n\,\sigma_0^{\frac{\la-1}{2}}\ \mid \ n\in \bZ\,\}\,,\ \ \text{and}
\\
&\qquad\qquad M_{\bC^*\!,\,\la,\text{odd}}\ \ \ \ \text{has basis}\ \ \ \{\,z^n\,\sigma_0^{\frac{\la-1}{2}}\ \mid \ n\in \bZ+ 1/2\,\}\,.
\end{aligned}
\end{equation}
For reasons of radial symmetry we only need to consider the integral of the expression (\ref{pol4}) over $\bC\cup \{\infty\}$ against the $U_\bR$-invariant measure $\omega$ when $f$ and $g$ are the same basis element. With the convention of (\ref{basis1}), with $n$ denoting either an integer or a true half integer, and using (\ref{URmetric1},\ref{URmetric2}), we find
\begin{equation}
\label{pol5}
\begin{aligned}
&(z^n\,\sigma_0^{\frac{\la-1}{2}}, z^n\,\sigma_0^{\frac{\la-1}{2}})_{\fu_\bR} \ = \ \int_{\bC\cup\{\infty\}}  |z|^{2n}\,\|\sigma_0\|^{\la-1}\,\omega \ =
\\
&\qquad\qquad\qquad = \ \int_{\bC\cup\{\infty\}} \frac{4\,|z|^{2n+\la-1}}{(1+ |z|^2)^{\la+1}}\,dz\,d\bar z\ = \ 8\,\pi\,\int_0^\infty\frac{r^{2n+\la}\,dr}{(1+ r^2)^{\la+1}}
\\
&\qquad\qquad\qquad\qquad\qquad\qquad = \ 4\,\pi\,\int_0^\infty\frac{u^{n+(\la-1)/2}\,du}{(1+ u)^{\la+1}}\ .
\end{aligned}
\end{equation}
This integral converges if and only if $-(\la + 1)/2 < n < (\la + 1)/2$. Since the integrand is positive, the integral has a strictly positive value in the range of convergence.

To continue the integral meromorphically beyond the range of convergence one uses integration by parts. Formally, for $s,t\in \bR$,
\begin{equation}
\label{pol6}
s\,\int_0^\infty \frac{u^{s-1}\,du}{(1+u)^t}\ \ = \ \ t\,\int_0^\infty \frac{u^s\,du}{(1+u)^{t+1}}\ .
\end{equation}
Applying this identity with $s=n+(\la+1)/2$ and $t=\la+1$, one finds that the integral changes sign and becomes strictly negative for $-(\la + 3)/2 < n < -(\la + 1)/2$. The same argument, with $-n$ substituted for $n$, shows that the integral is also strictly negative for $(\la+1)/2<n<(\la+3)/2$. Then the pattern continues: when $n$ is negative and decreased by one, or if $n$ is positive and increased by one, the integral changes sign. Poles occur when $\la$ is an odd integer in the even case, or an even integer in the odd case -- i.e., exactly when the module becomes reducible. That is what must happen, of course; the polarization is well defined only in the irreducible range. The preceding discussion can be summarized succinctly in terms of the Hodge filtration: with $\epsilon$ referring to either the even or the odd parity, in the irreducible range,
\begin{equation}
\label{pol7}
\begin{aligned}
&s \in F_0\,M_{\bC^*,\la,\epsilon}\ \ \Longrightarrow\ \ \text{the integral defining\,\ $(s,s)_{\fu_\bR}$\,\ converges}
\\
&\qquad\qquad s \in F_p\,M_{\bC^*,\la,\epsilon}\cap (F_p\,M_{\bC^*,\la,\epsilon})^\perp,\,\ s\neq 0\ \ \Longrightarrow\ \ (-1)^p\,(s,s)_{\fu_\bR}\, >\, 0\,;
\end{aligned}
\end{equation}
cf. (\ref{C*construction5}). The second statement is the assertion of our conjecture in the case of the open $K$-orbit $\bC^*$

The change in sign is directly related to a change in the weight filtration. Let $\la_0 > 0$ be reduction point -- i.e., a positive odd or even integer, depending on whether the parity is even or odd. In terms of the basis (\ref{basis1}), with the same convention of letting $n$ refer to an integer or true half integer, depending on the parity,
\begin{equation}
\label{basis2}
W_1\,M_{\bC^*\!,\,\la_0,\epsilon}\ \ \ \ \text{has basis}\ \ \ \{\,z^n\,\sigma_0^{\frac{\la-1}{2}}\ \mid \ -(\la - 1) \leq 2n \leq \la - 1\,\}\,.
\end{equation}
Thus, as the parameter $\la$ crosses the reduction point $\la_0$ going from left to right, the sign of $(z^n\,\sigma_0,z^n\,\sigma_0)_{\fu_\bR}$ remains the same if and only if $z^n\,\sigma_0 \in W_1\,M_{\bC^*\!,\,\la_0,\epsilon}$. This is one instance of a general fact. As the parameter for a family of induced representations crosses a reduction point, the sign changes of $(\ ,\ )$ are governed by the weight filtration -- that is the assertion of the Jantzen conjecture proved by Beilinson-Bernstein \cite{bb:1993}. The jumps of the Hodge filtration at the reduction point line up with the weight filtration, to produce exactly the sign changes predicted by our conjecture.

When $\la=m>0$ is a positive integer and $\epsilon$ denotes the opposite parity, i.e., the parity of $m+1$, the Harish Chandra module $W_1\,M_{\bC^*\!,\,m,\epsilon}$ has dimension $m$, and the integral (\ref{pol5}) converges for all basis elements: this is the usual description of the positive definite $U_\bR$-invariant inner product on the irreducible $m$-dimensional representation. For $\la=0$ and $\epsilon$ odd, $W_1\,M_{\bC^*\!,\,m,\epsilon}$ reduces to zero. The corresponding sheaf $W_1\,\cM_{\bC^*\!,\,0,\text{odd}}$ is the one and only irreducible Harish Chandra sheaf for $G_\bR = SU(1,1)$ without nonzero sections.

The two singleton orbits $\{0\}$, $\{\infty\}$ are related by an outer automorphism of $G_\bR = SU(1,1)$. Thus it is only necessary to discuss the orbit $\{0\}$. Since $K$ fixes the origin, it must act on the geometric fiber of any $\cD_\la$-module supported at $\{0\}$, and that forces an integral twisting parameter:
\begin{equation}
\label{pt1}
\la\ = \ m\,\in\, \bZ_{\geq 0}\,.
\end{equation}
In the untwisted case, the only irreducible $\cD$-module supported at the origin in $\bC$ is the one generated by the ``holomorphic delta function",
\begin{equation}
\label{pt2}
\cD_\bC\,\delta_0\ \ = \ \ \bC[z,z^{-1}]/\bC[z]\,.
\end{equation}
Thus $\de_0\cong z^{-1}$, and the $SL_2$-triple (\ref{sl2triple2}) acts according to the formulas
\begin{equation}
\label{pt3}
h\,\de_0\ = \ 2\de_0\,,\ \ \ e_-\,\de_0\ = \ 0\,,
\end{equation}
with $e_+$ acting freely, by normal differentiation. The section $\sigma_2$ of $\cL_2$ in (\ref{sl2triple2}) is nonzero at the origin and is $K$-invariant, and this leads to a description of the sheaf $\cM_{\{0\},m}$, or equivalently, to its space of global sections $M_{\{0\},m}$\,,
\begin{equation}
\label{pt4}
M_{\{0\},m}\ \ \ \ \text{has basis}\ \ \ \{\, ({\textstyle{\frac{d^n}{d z^n}}}\,\de_0)\,\sigma_2^{(m-1)/2}\ \mid \ n\geq 0\,\}\,;
\end{equation}
here $\sigma_2^{(m-1)/2}$ can be viewed as a section of $\cL_{m-1}$, a section that is regular and nonzero except at $\infty$. The $SL_2$-triple $e_+,\,h,\,e_-$ acts by the product rule, on $\frac{d^n}{d z^n}\,\de_0$ according to (\ref{sl2triple2}) and (\ref{pt3}), and on $\sigma_2^{(m-1)/2}$ according to (\ref{sl2triple3}). In particular,
\begin{equation}
\label{pt5}
\begin{aligned}
&e_+\left(({\textstyle{\frac{d^n}{d z^n}}}\,\de_0)\,\sigma_2^{(m-1)/2}\right)\ = \ -\,({\textstyle{\frac{d^{n+1}}{d z^{n+1}}}}\,\de_0)\,\sigma_2^{(m-1)/2}\,,
\\
&\qquad h\left(({\textstyle{\frac{d^n}{d z^n}}}\,\de_0)\,\sigma_2^{(m-1)/2}\right)\ = \ (2n+m+1)({\textstyle{\frac{d^{n+1}}{d z^{n+1}}}}\,\de_0)\,\sigma_2^{(m-1)/2}\,,
\\
&\qquad\qquad e_-\left(({\textstyle{\frac{d^n}{d z^n}}}\,\de_0)\,\sigma_2^{(m-1)/2}\right)\ = \ n(n+1)({\textstyle{\frac{d^{n-1}}{d z^{n-1}}}}\,\de_0)\,\sigma_2^{(m-1)/2}\,,
\end{aligned}
\end{equation}
as can be checked readily.

The inclusion $\{0\} \hookrightarrow \bC \cup \{\infty\}$ is a very special case of a closed embedding. In general the $\cD$-module direct image of an irreducible module under a closed embedding remains irreducible, so the weight filtration collapses. The effect of closed embeddings on the Hodge filtration also has a simple description: the Hodge index is increased by the order of normal derivative. In the case of $M_{\{0\},m}$ this means
\begin{equation}
\label{pt6}
\begin{aligned}
&W_0\,M_{\{0\},m}\ \ = \ \ M_{\{0\},m}\ ,
\\
&\qquad F_p\,M_{\{0\},m}\ \ \ \ \text{has basis}\ \ \ \{\, ({\textstyle{\frac{d^n}{d z^n}}}\,\de_0)\,\sigma_2^{(m-1)/2}\ \mid \ 0\leq n\leq p-1\,\}\,,
\end{aligned}
\end{equation}
because the weight equals the dimension of the support, and the Hodge filtration starts at the codimension of the support.

The polarization pairs $\de_0$ and $\overline{\de_0}$ into $\de_{\bR,0}$, the delta function in the usual sense on $\bC\cong \bR^2$. It also pairs $\sigma_2^{(m-1)/2}$ and its complex conjugate into $\|\sigma_2^{(m-1)/2}\|^{2}= \|\sigma_2\|^{m-1} = (1+|z|^2)^{-m+1}$, as follows from (\ref{URmetric2}). Thus
\begin{equation}
\label{pt7}
\begin{aligned}
&(\,({\textstyle{\frac{d^k}{d z^k}}}\,\de_0)\sigma_2^{(m-1)/2}\,,\,\overline{({\textstyle{\frac{d^\ell}{d z^\ell}}}\,\de_0)\sigma_2^{(m-1)/2}}\,)_{\fu_\bR}= \int_{\bC\cup\{\infty\}}\!{\textstyle{\frac{d^k}{d z^k}}}\,{\textstyle{\frac{d^\ell}{d \bar z^\ell}}}\,\|\sigma_2^{(m-1)/2}\|^{2}\,\de_{\bR,0}\,\,\om \!\!
\\
&\ \, = \int_{\bC\cup\{\infty\}}\!{\textstyle{\frac{d^k}{d z^k}}}\,{\textstyle{\frac{d^\ell}{d \bar z^\ell}}}\,(1+|z|^2)^{-m-1}\,\de_{\bR,0}\,\,dz\,d\bar z\,\ = \ \left.{\textstyle{\frac{d^k}{d z^k}}}\,{\textstyle{\frac{d^\ell}{d \bar z^\ell}}}\,(1+|z|^2)^{-m-1}\right|_{z=0} \!\!
\end{aligned}
\end{equation}
vanishes unless $k=\ell$, in which case
\begin{equation}
\label{pt8}
(\,({\textstyle{\frac{d^k}{d z^k}}}\,\de_0)\sigma_2^{(m-1)/2}\,,\,\overline{({\textstyle{\frac{d^k}{d z^k}}}\,\de_0)\sigma_2^{(m-1)/2}}\,)_{\fu_\bR}\ =\ (-1)^k\,k!\ {\prod}_{j=1}^k\,(m+j)\ .
\end{equation}
That, of course, is consistent with our conjecture in this particular instance.

\end{document}